\documentclass[12pt]{article}

\usepackage{amsmath,amssymb}
\usepackage{tikz}
\usetikzlibrary{decorations.pathreplacing}

\usetikzlibrary{calc}
\usepackage{color}

\newtheorem{theorem}{Theorem}[section]
\newtheorem{lemma}[theorem]{Lemma}

\newtheorem{proposition}[theorem]{Proposition}
\newtheorem{fact}[theorem]{Fact}

\newtheorem{remark}[theorem]{Remark}

\tikzstyle{noeud}=[circle,inner sep=2, minimum size =3 pt, line width = 1pt, draw=black, fill=white]

\newcommand{\proof}{\noindent{\bf Proof.\ }}
\newcommand{\qed}{\hfill $\square$ \bigskip}
\newcommand{\Dom}{{Dominator}}
\newcommand{\St}{{Staller}}
\newcommand{\gmb}{\gamma_{{\rm MB}}}
\def\cp{\,\square\,}

\newcommand\litem[1]{\item{\textit{Case #1:}}}

\begin{document}

\title{Maker-Breaker domination number}

\author{
    Valentin Gledel $^{a}$
    \and
	Vesna Ir\v si\v c $^{b,c}$
	\and
	Sandi Klav\v zar $^{b,c,d}$
}

\date{\today}

\maketitle
\begin{center}
	$^a$ Univ Lyon, Universit\'e Lyon 1, LIRIS UMR CNRS 5205, F-69621, Lyon, France \\
	\medskip

	$^b$ Institute of Mathematics, Physics and Mechanics, Ljubljana, Slovenia\\
	\medskip

	$^c$ Faculty of Mathematics and Physics, University of Ljubljana, Slovenia\\
	\medskip
	
	$^d$ Faculty of Natural Sciences and Mathematics, University of Maribor, Slovenia
		
\end{center}

\begin{abstract}
The Maker-Breaker domination game is played on a graph $G$ by \Dom\ and \St. The players alternatively select a vertex of $G$ that was not yet chosen in the course of the game. \Dom\ wins if at some point the vertices he has chosen form a dominating set. \St\ wins if \Dom\ cannot form a dominating set. In this paper we introduce the Maker-Breaker domination number $\gmb(G)$ of $G$ as the minimum number of moves of Dominator to win the game provided that he has a winning strategy and is the first to play. If Staller plays first, then the corresponding invariant is denoted $\gmb'(G)$. Comparing the two invariants it turns out that they behave much differently than the related game domination numbers. The invariant $\gmb(G)$ is also compared with the domination number. Using the Erd\H os-Selfridge Criterion a large class of graphs $G$ is found for which $\gmb(G) > \gamma(G)$ holds. Residual graphs are introduced and used to bound/determine $\gmb(G)$ and $\gmb'(G)$. Using residual graphs, $\gmb(T)$ and $\gmb'(T)$ are determined for an arbitrary tree. The invariants are also obtained for cycles and bounded for union of graphs. A list of open problems and directions for further investigations is given. 
\end{abstract}


\noindent {\bf Key words:} Maker-Breaker domination game; Maker-Breaker domination number; domination game; perfect matching; tree; cycle; union of graphs  

\medskip\noindent
{\bf AMS Subj.\ Class:} 05C57, 05C69, 91A43

\section{Introduction}
\label{sec:intro}

{\em Maker-Breaker games} (as well as other {\em positional games}) have been introduced by Erd\H os and Selfridge in~\cite{erdos-1973}, and since then have been the subject of numerous studies, see~\cite{bednarska-2000, ben-shimon-2011, ferber-2017, hefetz-2008}. Maker-Breaker games are played on hypergraphs by two players called Maker and Breaker. They take turns and at each turn the current player selects a new vertex. Maker wins if at some point of the game he has selected all vertices from one of the hyperedges, while Breaker wins if she can keep him from doing it. See~\cite{beck-2008} and~\cite{hefetz-2014} for general introductions on this field. 

Very recently, the {\em Maker-Breaker domination game} was introduced in~\cite{gledel-2018+}. The game is played on a graph $G$ with two players named \Dom\ and \St. These names were selected to emphasize the domination nature of the game and to be consistent with the usual domination game where these two names are standard by now. (The domination game was introduced in~\cite{brklra-2010} and further studied in dozens of papers, cf.~\cite{bu-2015, dorbec-2015, NSS2016, sc-2016, XLK-2018}.) The players alternatively select a vertex of $G$ that was not yet chosen in the course of the game. \Dom\ wins if at some point, the vertices he has chosen form a dominating set. \St\ wins if \Dom\ cannot form a dominating set. 
Note that the Maker-Breaker domination game is a Maker-Breaker game. Indeed, if for a graph $G$ we build a hypergraph $\mathcal F$ with the same set of vertices as $G$, and in which the hyperedges are the dominating sets of $G$, then Dominator wins the Maker-Breaker domination game on $G$ if and only if Maker wins the Maker-Breaker game on $\mathcal F$. 

In several papers on Maker-Breaker games the authors were interested in the smallest number of moves needed for Maker to win, see~\cite{clemens-2012, clemens-2018, hefetz-2008}. Also, in~\cite{gledel-2018+} it was emphasized that when dealing with the Maker-Breaker games, there are two natural questions: (i) which player has a winning strategy and (ii) what is the minimum number of moves if Dominator has a winning strategy. In the seminal paper question (i) is investigated, while in this paper we study (ii). For this sake we say that if $G$ is a graph, then the {\em Maker-Breaker domination number} $\gmb(G)$ of $G$ is the minimum number of moves of Dominator to win the game provided that he has a winning strategy and is the first to play. Otherwise we set $\gmb(G) = \infty$. Similarly, $\gmb'(G)$ denotes is the minimum number of moves of Dominator in the game in which Staller plays first. 

We proceed as follows. In the next section we list additional definitions and several known results needed in this paper, as well as prove some basic results on the Maker-Breaker domination number. In Section~\ref{sec:MB-numbers} we first compare $\gmb(G)$ with $\gmb'(G)$ and find out that they behave totally different than the related game domination invariants. We also compare  $\gmb(G)$ with the domination number and using the Erd\H os-Selfridge Criterion prove that if the number of $\gamma$-sets of $G$ is not too big, then $\gmb(G) > \gamma(G)$. In Section~\ref{sec:residual} we introduce residual graphs, determine (resp.\ bound) $\gmb'(G)$ (resp.\ $\gmb(G)$) in terms of the residual graph, and determine $\gmb(T)$ and $\gmb'(T)$ for an arbitrary tree. In the next two sections we obtain the invariants for cycles and bound them for union of graphs. We conclude with a list of open problems and directions for further investigation.

\section{Preliminaries}
\label{sec:prelim}

Let $G$ be a graph. A vertex of $G$ adjacent to a leaf is a {\em support vertex} of $G$. A {\em perfect matching} of $G$ is a set of pairwise independent edges that cover $V(G)$. The order of $G$ will be denoted with $n(G)$. If $u$ is a vertex of $G$, then $N[u]$ denotes the {\em closed neighborhood} of $u$. If $v$ is another vertex then we set $N[u,v] = N[u]\cap N[v]$. A set $D\subseteq V(G)$ is a {\em dominating set} of $G$ if $\cup_{u\in D}N[u] = V(G)$. The {\em domination number} $\gamma(G)$ is the size of a smallest dominating set of $G$. A dominating set of size $\gamma(G)$ is called a {\em $\gamma$-set} of $G$. 

The Maker-Breaker domination game is called a D-game (resp.\ S-game) if \Dom\ (resp.\  \St) is the first to play a vertex. The sequence of vertices selected in a D-game will be denoted with $d_1, s_1, d_2, s_2, \ldots$, and the sequence of vertices selected in an S-game with $s_1', d_1', s_2', d_2', \ldots$
Suppose that \Dom\ wins a D-game. Then the last vertex played is by \Dom, let it be $d_k$. By the definition of the game, $\{d_1, \ldots, d_k\}$ is a dominating set of $G$.  Similarly, if \Dom\ wins an S-game and the last vertex played by \Dom\ is $d_\ell'$, then $\{d_1', \ldots, d_\ell'\}$ is a dominating set of $G$. 

Let $G$ be a graph, $k\ge 1$, and $u_1, \ldots, u_k, v_1,\ldots, v_k$ pairwise different vertices of $G$. Then we say that $X = \{ \{u_1,v_1\},\ldots, \{u_k,v_k\} \}$ is a {\em pairing dominating set} if 
$$\bigcup_{i=1}^k N[u_i,v_i] = V(G)\,.$$
In the rest we will use this concept via the following interpretation proved in~\cite[Proposition 9]{gledel-2018+}. To be self-contained, we give here an alternative, short proof.

\begin{lemma}
\label{lem:dominating-pairs}
Let $u_1, \ldots, u_k, v_1,\ldots, v_k$ be pairwise different vertices of a graph $G$, and let $X = \{ \{u_1,v_1\},\ldots, \{u_k,v_k\} \}$. Then $X$ is a pairing dominating set if and only if every set $\{x_1,\ldots, x_k\}$, where $x_i\in \{u_i, v_i\}$, $i\in [k]$, is a dominating set of $G$.   
\end{lemma}

\proof
Suppose first that $X$ is a pairing dominating set, that is, $\cup_{i=1}^k N[u_i,v_i] = V(G)$. Let $\{x_1,\ldots, x_k\}$ be an arbitrary set with $x_i\in \{u_i, v_i\}$, $i\in [k]$. Then $V(G) = \cup_{i=1}^k N[u_i,v_i] \subseteq \cup_{i=1}^k N[x_i]$. So $\{x_1,\ldots, x_k\}$ is a dominating set of $G$. 

Conversely, consider a set $\{x_1,\ldots, x_k\}$, where $x_i\in \{u_i, v_i\}$, $i\in [k]$, and suppose that $\cup_{i=1}^k N[u_i,v_i]$ is a proper subset of $V(G)$. Let $w\in V(G)\setminus \cup_{i=1}^k N[u_i,v_i]$. Then for every $i\in [k]$ we have $w\notin N[u_i,v_i]$. Let $y_i\in \{u_i,v_i\}$ be such that $w\notin N[y_i]$. But then $\cup_{i=1}^kN[y_i]$ is a proper subset of $V(G)$, that is, $\{y_1,\ldots, y_k\}$ is not a dominating set.  
\qed

If $X = \{ \{u_1,v_1\},\ldots, \{u_k,v_k\} \}$ is a pairing dominating set such that $u_iv_i\in E(G)$ holds for $i\in [k]$, then we say that $X$ is a {\em dominating matching}. 

\begin{fact}  {\rm \cite[Proposition 10]{gledel-2018+}}
\label{fact:dominating-pair}
If $G$ admits a pairing dominating set, then Dominator has a winning strategy on $G$ in the D-game as well as in the S-game. 
\end{fact}

The converse of Fact~\ref{fact:dominating-pair} does not hold in general. For instance, in~\cite[Figure 4]{gledel-2018+} a chordal graph is presented on which Dominator has a winning strategy in both games but which admits no pairing dominating set. On the other hand, the converse holds in the class of trees because if Dominator has a winning strategy on a tree $T$, then it was proved in~\cite{gledel-2018+} that $T$ has a dominating matching. Moreover, the converse also holds for co-graphs. The next lemma considers the variation of the game where players might skip some moves. This means that the current player selects no vertex and the previous player just makes another move. The skipped moves do not count in $\gmb(G)$ or $\gmb'(G)$.

\begin{lemma} {\rm (No-Skip Lemma)}
\label{lem:no-skip}
In an optimal strategy of Dominator to achieve $\gmb(G)$ or $\gmb'(G)$ it is never an advantage for him to skip a move. Moreover, if Staller skips a move it can never disadvantage Dominator. 
\end{lemma}

\proof
Suppose a D-game or an S-game is played. Let Dominator and Staller play optimally until some point when Staller decides to skip a move. In that case, Dominator imagines an arbitrary move of Staller, say $x$, and replies optimally to this move. Since Dominator can always, no matter the way Staller selects vertices, finish the game in no more that $\gmb(G)$ (resp.\ $\gmb'(G)$) moves, this property is preserved after the imagined move $x$ and the reply to it. Then Dominator proceeds until the end of the game with the same strategy. Note that it may happen that in the course of the game Staller selects a vertex which is not a legal move in the game Dominator is imagining. In that case Dominator imagines that yet some other legal move has been played by Staller. In this way the game on $G$ will finish in no more than $\gmb(G)$ (resp.\ $\gmb'(G)$) moves. 

With a strategy of Staller parallel to the above strategy of Dominator we also infer that it is never an advantage for Dominator to skip a move. 
\qed

If $G$ is a graph and $S\subseteq V(G)$, then let $G|S$ denote that graph $G$ in which the vertices from $S$ are declared to be already dominated, that is, Dominator is not obliged to dominate them in the rest of the game. Then we have the following Continuation Principle, a proof of which is much simpler that the corresponding principle for the domination game~\cite{kinnersley-2013}. 

\begin{remark} {\rm (Continuation Principle)}
\label{rem:continuation-principle}
Let $G$ be a graph with $A,B\subseteq V(G)$. If $B\subseteq A$, then $\gmb(G|A)\leq \gmb(G|B)$ and $\gmb'(G|A)\leq \gmb'(G|B)$.
\end{remark}

Indeed, the remark follows from the fact that Dominator can apply the same strategy in $G|A$ as in $G|B$. 

Suppose that $\gmb(G) < \infty$. Then in any winning strategy of Dominator, he will play at most half of the vertices (because Staller will play the other half) which in turn implies that 
\begin{equation}
\label{eq:n/2-Dominator}
1\le \gmb(G)\le \left\lceil \frac{n(G)}{2}\right\rceil\,.
\end{equation}
The bound is sharp, consider for instance the disjoint union of $K_1$ and several copies of $K_2$. It is also easy to see that all the possible values from~\eqref{eq:n/2-Dominator} can be realized by considering the disjoint union of a complete graph and an appropriate number of $K_2$s. Similarly, for the S-game, assuming that $\gmb'(G) < \infty$,  we have
\begin{equation}
\label{eq:n/2-Staller}
1\le \gmb'(G)\le \left\lfloor \frac{n(G)}{2}\right\rfloor\,,
\end{equation}
where again all the values can be realized. 

Later we will apply the celebrated Erd\H os-Selfridge Criterion for Maker-Breaker games that reads as follows.

\begin{theorem}[Erd\H os-Selfridge Criterion~\cite{erdos-1973}]
\label{thm:criterion}
If $\mathcal F$ is a hypergraph, then
$$ \sum_{A \in \mathcal F} 2^{-|A|}<\frac{1}{2} \ \Rightarrow \ \text{$\mathcal{F}$ is a Breaker's win}\,.$$
\end{theorem}

This theorem together with its proof can also be found in the book~\cite[Theorem 2.3.3]{hefetz-2014}.

\section{Maker-Breaker domination numbers}
\label{sec:MB-numbers}

In this section we first compare $\gmb(G)$ with $\gmb'(G)$ and construct graphs for all possible values of the invariants. In the second part we compare $\gmb(G)$ with the domination number and using the Erd\H os-Selfridge Criterion find a large class of graphs $G$ for which $\gmb(G) > \gamma(G)$ holds.

\subsection{Realizations of Maker-Breaker domination numbers}

One of the fundamental theorems on the domination game proved in~\cite{brklra-2010, kinnersley-2013} asserts that $|\gamma_g(G) - \gamma_g'(G)|\le 1$ holds for every graph $G$. The next result reveals that the situation with the Maker-Breaker domination number is dramatically different.   

\begin{theorem}
\label{thm:realization}
If $G$ is a graph, then $\gamma(G) \leq \gmb(G) \leq \gmb'(G)$. Moreover, for any integers $r,s,t$, where $2\le r\le s\le t$, there exists a graph $G$ such that $\gamma(G) = r$, $\gmb(G) = s$, and $\gmb'(G) = t$.  
\end{theorem}

\proof
The assertion $\gamma(G) \leq \gmb(G)$ is clear since $\{d_1, d_2, \ldots\}$ is a dominating set of $G$. 

A D-game can be viewed as an S-game in which Staller has skipped her first move. The No-Skip Lemma thus implies that $\gmb(G) \leq \gmb'(G)$. 

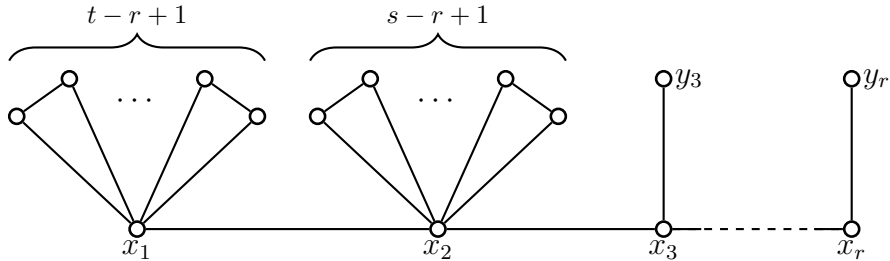
\begin{figure}[ht!]
\begin{center}
\begin{tikzpicture}[scale=1.0,style=thick,x=1cm,y=1cm]
\def\vr{3pt}


\node[noeud] (x1) at (0,0){};
\node[noeud] (x2) at (4,0){};
\node[noeud] (x3) at (7,0){};
\node[noeud] (xr) at (9.5,0){};
\node[noeud] (x3') at (7,2){};
\node[noeud] (xr') at (9.5,2){};
\node[noeud] (y1) at (-1.6,1.5){};
\node[noeud] (y1') at (-0.9,2){};
\node[noeud] (yt-k-1) at (0.9,2){};
\node[noeud] (yt-k-1') at (1.6,1.5){};
\node[noeud] (z1) at (2.4,1.5){};
\node[noeud] (z1') at (3.1,2){};
\node[noeud] (zt-k-1) at (4.9,2){};
\node[noeud] (zt-k-1') at (5.6,1.5){};

\draw (x1) -- (x2) -- (x3);
\draw (x1) -- (y1) -- (y1') -- (x1);
\draw (x1) -- (yt-k-1) -- (yt-k-1') -- (x1);
\draw (x2) -- (z1) -- (z1') -- (x2);
\draw (x2) -- (zt-k-1) -- (zt-k-1') -- (x2);
\draw (x3) -- (7.5,0);
\draw[dashed] (x3) -- (xr);
\draw (xr) -- (9,0);
\draw (x3) -- (x3');
\draw (xr) -- (xr');


\draw [decorate,decoration={brace,amplitude=10pt,raise=4pt}] (-1.7,2.1)-- (1.7,2.1)node[above= 13pt,midway]{\footnotesize $t-r+1$};
\draw [decorate,decoration={brace,amplitude=10pt,raise=4pt}] (2.3,2.1)-- (5.7,2.1)node[above= 13pt,midway]{\footnotesize $s-r+1$};
\draw (0,1.7) node {$\cdots$};
\draw (4,1.7) node {$\cdots$};
\draw(x1) node[below] {$x_{1}$};
\draw(x2) node[below] {$x_{2}$};
\draw(x3) node[below] {$x_{3}$};
\draw(xr) node[below] {$x_{r}$};
\draw(x3') node[right] {$y_{3}$};
\draw(xr') node[right] {$y_{r}$};

\end{tikzpicture}
\end{center}
\caption{Graph $G_{r,s,t}$}
\label{fig:Grst}
\end{figure}

Let $r, s, t$ be fixed integers where $2\le r\le s\le t$. Construct a graph $G_{r,s,t}$ as follows. Start with a path of length $r-1$ on consecutive vertices $x_1,\ldots, x_r$. Attach $t-r+1$ pendant triangles at $x_1$ and $s-r+1$ pendant triangles at $x_2$. Finally, at each (if any) of the vertices $x_3,\ldots, x_r$ attach a pendant vertex $y_3,\ldots, y_r$, respectively. The construction should be clear with the aid of Fig.~\ref{fig:Grst}. 

The $t-r+1$ triangle edges opposite to $x_1$ and the $s-r+1$ triangle edges opposite to $x_2$ together with the edges $x_3y_3,\ldots, x_ry_r$ form a dominating matching of 
$G_{r,s,t}$. (If one would like to have a matching, then adding the edge $x_1x_2$ would do the job.) Hence by Fact~\ref{fact:dominating-pair}, Dominator has a winning strategy in both games. We claim that $\gamma(G_{r,s,t}) = r$, $\gmb(G_{r,s,t}) = s$, and $\gmb'(G_{r,s,t}) = t$, where the first assertion is clear. 

Consider the D-game. We first describe the following strategy of Dominator. He starts the game with the move $d_1 = x_1$. Then no matter how Staller plays, Dominator can proceed such that at most $(r-2) + (s-r+1)$ moves of him will be needed to dominate the graph. To do this, whenever Staller plays on $x_i$ or $y_i$, $i\ge 3$, Dominator replies with $y_i$ or $x_i$, respectively. Also, Dominator proceeds along the same lines when Staller plays on a vertex that lies in a triangle attached to $x_2$. Using this strategy Dominator ensures that he will play at most $1 + (r-2) + (s-r+1) = s$ moves which means that $\gmb(G) \le s$. To prove the other inequality, consider the following strategy of Staller. If Dominator starts with $d_1 = x_1$, then Staller sets $s_1 = x_2$. Then Dominator will need at least $(r-2) + (s-r+1)$  additional moves to dominate $G$. On the other hand, if $d_1 \ne x_1$, then Staller sets $s_1 = x_1$, but then Dominator will need at least $(r-2) + (t-r+1) \ge (r-2) + (s-r+1)$ additional moves. In any case,  $\gmb(G) \ge  1 + (r-2) + (s-r+1) = s$.

Consider next the S-game. The proof that $\gmb'(G_{r,s,t}) = t$ proceeds similarly as above. To show that $\gmb(G) \ge t$, Staller can apply a strategy to start the S-game with $s_1' = x_1$. In this way she can guarantee that Dominator will need to play at least $(t-r+1) + 1 + (r-2) = t$ moves. On the other hand, Dominator can play such that no more than $t$ vertices will be selected by him. Whenever Staller plays $x_i$ or $y_i$, $i\ge 3$, he replies with $y_i$ or $x_i$, respectively. Moreover, during the game he will be able to play $x_1$ or $x_2$.  
\qed

Note that if $\gamma(G) = 1$, then also $\gmb(G) = 1$. Hence Theorem~\ref{thm:realization} does not extend to the case $r=1$. On the other hand, if $G_t$, $t\ge 1$, is the graph obtained from $t$ disjoint triangles by identifying a vertex from each of the triangles (so that this new vertex is of degree $2t$), then $\gamma(G_t) = 1$, $\gmb(G_t) = 1$, and $\gmb'(G_t) = t$. 

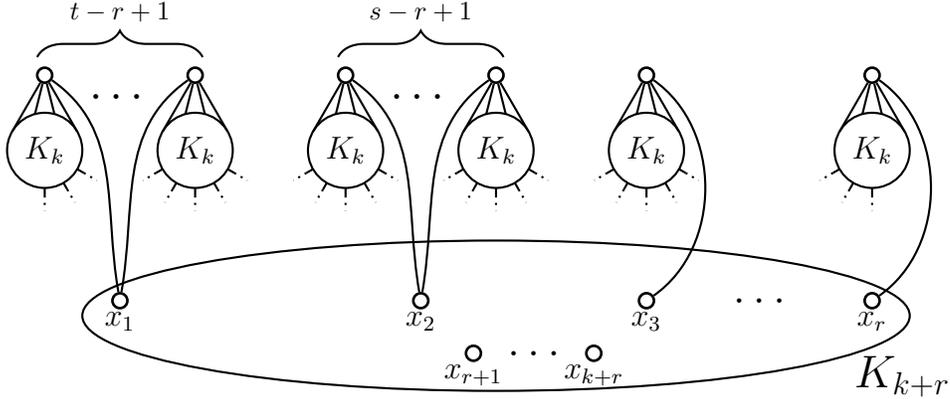
\begin{figure}[ht!]
\begin{center}
\begin{tikzpicture}[scale=1.0,style=thick,x=1cm,y=1cm]

\node[noeud] (x1) at (0,0){};
\node[noeud] (x2) at (4,0){};
\node[noeud] (x3) at (7,0){};
\node[noeud] (xr) at (10,0){};

\node[noeud] (xr+1) at (4.7,-0.7){};
\node[noeud] (xk+r) at (6.3,-0.7){};
\node[below] at (xr+1){$x_{r+1}$};
\node[below] at (xk+r){$x_{k+r}$};

\draw (-1,2) circle (0.5);
\node at (-1,2) {$K_k$};
\node[noeud] (y11) at (-1,3){};
\draw (y11) to[out=-35, in =100] (x1);
\begin{scope}[xshift=-1cm, yshift=2cm]
\draw (270:0.5) -- (270:0.6);
\draw (300:0.5) -- (300:0.6);
\draw (330:0.5) -- (330:0.6);
\draw[dotted] (270:0.6) -- (270:0.8);
\draw[dotted] (300:0.6) -- (300:0.8);
\draw[dotted] (330:0.6) -- (330:0.8);
\draw (30:0.5) -- (y11);
\draw (75:0.5) -- (y11);
\draw (105:0.5) -- (y11);
\draw (150:0.5) -- (y11);
\end{scope}

\draw (1,2) circle (0.5);
\node at (1,2) {$K_k$};
\node[noeud] (y12) at (1,3){};
\draw (y12) to[out=215, in =80] (x1);
\begin{scope}[xshift=1cm, yshift=2cm]
\draw (210:0.5) -- (210:0.6);
\draw (240:0.5) -- (240:0.6);
\draw (270:0.5) -- (270:0.6);
\draw (300:0.5) -- (300:0.6);
\draw (330:0.5) -- (330:0.6);
\draw[dotted] (210:0.6) -- (210:0.8);
\draw[dotted] (240:0.6) -- (240:0.8);
\draw[dotted] (270:0.6) -- (270:0.8);
\draw[dotted] (300:0.6) -- (300:0.8);
\draw[dotted] (330:0.6) -- (330:0.8);
\draw (30:0.5) -- (y12);
\draw (75:0.5) -- (y12);
\draw (105:0.5) -- (y12);
\draw (150:0.5) -- (y12);
\end{scope}

\draw (3,2) circle (0.5);
\node at (3,2) {$K_k$};
\node[noeud] (y21) at (3,3){};
\draw (y21) to[out=-35, in =100] (x2);
\begin{scope}[xshift=3cm, yshift=2cm]
\draw (210:0.5) -- (210:0.6);
\draw (240:0.5) -- (240:0.6);
\draw (270:0.5) -- (270:0.6);
\draw (300:0.5) -- (300:0.6);
\draw (330:0.5) -- (330:0.6);
\draw[dotted] (210:0.6) -- (210:0.8);
\draw[dotted] (240:0.6) -- (240:0.8);
\draw[dotted] (270:0.6) -- (270:0.8);
\draw[dotted] (300:0.6) -- (300:0.8);
\draw[dotted] (330:0.6) -- (330:0.8);
\draw (30:0.5) -- (y21);
\draw (75:0.5) -- (y21);
\draw (105:0.5) -- (y21);
\draw (150:0.5) -- (y21);
\end{scope}

\draw (5,2) circle (0.5);
\node at (5,2) {$K_k$};
\node[noeud] (y22) at (5,3){};
\draw (y22) to[out=215, in =80] (x2);
\begin{scope}[xshift=5cm, yshift=2cm]
\draw (210:0.5) -- (210:0.6);
\draw (240:0.5) -- (240:0.6);
\draw (270:0.5) -- (270:0.6);
\draw (300:0.5) -- (300:0.6);
\draw (330:0.5) -- (330:0.6);
\draw[dotted] (210:0.6) -- (210:0.8);
\draw[dotted] (240:0.6) -- (240:0.8);
\draw[dotted] (270:0.6) -- (270:0.8);
\draw[dotted] (300:0.6) -- (300:0.8);
\draw[dotted] (330:0.6) -- (330:0.8);
\draw (30:0.5) -- (y22);
\draw (75:0.5) -- (y22);
\draw (105:0.5) -- (y22);
\draw (150:0.5) -- (y22);
\end{scope}

\draw (7,2) circle (0.5);
\node at (7,2) {$K_k$};
\node[noeud] (y3) at (7,3){};
\draw (y3) to[out=-35, in =35] (x3);
\begin{scope}[xshift=7cm, yshift=2cm]
\draw (210:0.5) -- (210:0.6);
\draw (240:0.5) -- (240:0.6);
\draw (270:0.5) -- (270:0.6);
\draw (300:0.5) -- (300:0.6);
\draw (330:0.5) -- (330:0.6);
\draw[dotted] (210:0.6) -- (210:0.8);
\draw[dotted] (240:0.6) -- (240:0.8);
\draw[dotted] (270:0.6) -- (270:0.8);
\draw[dotted] (300:0.6) -- (300:0.8);
\draw[dotted] (330:0.6) -- (330:0.8);
\draw (30:0.5) -- (y3);
\draw (75:0.5) -- (y3);
\draw (105:0.5) -- (y3);
\draw (150:0.5) -- (y3);
\end{scope}

\draw (10,2) circle (0.5);
\node at (10,2) {$K_k$};
\node[noeud] (yr) at (10,3){};
\draw (yr) to[out=-35, in =35] (xr);
\begin{scope}[xshift=10cm, yshift=2cm]
\draw (210:0.5) -- (210:0.6);
\draw (240:0.5) -- (240:0.6);
\draw (270:0.5) -- (270:0.6);
\draw[dotted] (210:0.6) -- (210:0.8);
\draw[dotted] (240:0.6) -- (240:0.8);
\draw[dotted] (270:0.6) -- (270:0.8);
\draw (30:0.5) -- (yr);
\draw (75:0.5) -- (yr);
\draw (105:0.5) -- (yr);
\draw (150:0.5) -- (yr);
\end{scope}

\node at (10.4,-1) {\Large $K_{k+r}$};

\draw [decorate,decoration={brace,amplitude=10pt,raise=4pt}] (-1.1,3.1)-- (1.1,3.1)node[above= 13pt,midway]{\footnotesize $t-r+1$};
\draw [decorate,decoration={brace,amplitude=10pt,raise=4pt}] (2.9,3.1)-- (5.1,3.1)node[above= 13pt,midway]{\footnotesize $s-r+1$};
\draw (0,2.7) node {\Large $\cdots$};
\draw (4,2.7) node {\Large $\cdots$};
\draw (8.5,0) node {\Large $\ldots$};
\draw (5.5,-0.7) node {\Large $\ldots$};
\draw(x1) node[below] {$x_{1}$};
\draw(x2) node[below] {$x_{2}$};
\draw(x3) node[below] {$x_{3}$};
\draw(xr) node[below] {$x_{r}$};

\draw (5,-0.2) ellipse (5.5 and 1);

\end{tikzpicture}
\end{center}
\caption{Graph $H_{k,r,s,t}$}
\label{fig:Hkrst}
\end{figure}

Theorem~\ref{thm:realization} extends also to highly connected graphs. To see this, consider the graphs $H_{k,r,s,t}$, $2\le r\le s\le t$, $k\ge 1$, that are schematically drawn in Fig.~\ref{fig:Hkrst}. Here, each vertex of a $K_k$ clique is adjacent to each vertex of the clique $K_{k+r}$. Then by arguments similar to those from the proof of Theorem~\ref{thm:realization} one can see that $\gamma(H_{k,r,s,t}) = r$, $\gmb(H_{k,r,s,t}) = s$, and $\gmb'(H_{k,r,s,t}) = t$. Moreover, $H_{k,r,s,t}$ is $(k+1)$-connected. 

\subsection{Relation with the domination number}
\label{sec:MB-gamma}

As already observed above, $\gamma(G) = 1$ if and only if $\gmb(G) = 1$. In general it would be interesting to characterize the graphs $G$ such that $\gmb(G) = \gamma(G) = k$, where $k\ge 2$ is a fixed integer. For $k=2$ the answer is simple: 

\begin{proposition}
\label{prp:gmb=gamma=2}
Let $G$ be a graph with $\gamma(G) = 2$. Then $\gmb(G) = \gamma(G) = 2$ if and only if $G$ has a vertex that lies in at least two $\gamma$-sets of $G$.
\end{proposition}

\proof
Suppose $\gmb(G) = \gamma(G) = 2$. After the moves $d_1$ and $s_1$ are played, Dominator has a strategy to win the game with the move $d_2$. Then $\{d_1, d_2\}$ is a $\gamma$-set of $G$. Moreover, Staller has an option to select $s_1$ such that $\{d_1,s_1\}$ is a $\gamma$-set, hence $d_1$ must lie in at least two $\gamma$-sets. 

Conversely, let $u$ be a vertex that lies in two $\gamma$-sets of $G$. Then Dominator plays $d_1 = u$, and then no matter which vertex is selected by Staller in her first move, Dominator can finish the game in his second move. 
\qed

Proposition~\ref{prp:gmb=gamma=2} can be rephrased to hold for larger $k$ also, but this would be more or less just rephrasing the definitions. It would be more interesting to find a structural characterization of the corresponding graphs. This task, however, seems difficult. On the other hand, the Erd\H os-Selfridge Criterion gives a sufficient condition for $\gmb(G) > \gamma(G)$. Let $X_\gamma(G)$ be the number of $\gamma$-sets of a graph $G$, cf.~\cite{connolly-2016}. Then: 

\begin{proposition}
\label{prop:gmb>gamma}
If $G$ is a graph and $X_\gamma(G) < 2^{\gamma(G)-1}$, then $\gmb(G) > \gamma(G)$.
\end{proposition}

\proof
Let $\mathcal F$ be the hypergraph with $V({\mathcal F}) = V(G)$ and whose hyperedges are the $\gamma$-sets of $G$. Then Theorem~\ref{thm:criterion} asserts that
\begin{equation}
\label{eq:erdos-for-gamma-sets}
\sum_{A \in \mathcal F} 2^{-|A|}<\frac{1}{2} \ \Rightarrow \ \text{$\mathcal{F}$ is a Breaker's win}\,.
\end{equation}
Since $|E(\mathcal F)| = X_\gamma(G)$ and each of these hyperedges has size $\gamma(G)$, we can estimate as follows: 
$$\sum_{A \in \mathcal F} 2^{-|A|} = \sum_{A \in \mathcal F} 2^{-\gamma(G)} = X_\gamma(G)\cdot 2^{-\gamma(G)} < 2^{\gamma(G)-1}\cdot 2^{-\gamma(G)} = \frac{1}{2}\,.$$
Therefore, $\mathcal{F}$ is a Breaker's win by~\eqref{eq:erdos-for-gamma-sets}. But this means that in $G$, Dominator is unable to win with $\gamma(G)$ moves, thus $\gmb(G) > \gamma(G)$.
\qed

Consider the cycles $C_{3k - 1}$, $k \geq 1$. It is known and easy to see that $\gamma(C_{3k-1}) = k$. We now determine the number of $\gamma$-sets of $C_{3k-1}$. Each vertex from a $\gamma$-set dominates itself and its two neighbors. As there are $k$ such triplets and $3k-1$ vertices in the graph, there is only one vertex that is dominated by two vertices from the $\gamma$-set, all others are dominated exactly once. Thus if the vertex that is dominated twice is fixed, then the $\gamma$-set of the cycle is uniquely determined. As there are $3k-1$ choices for this vertex, we have $X_\gamma(C_{3k - 1}) = 3k-1$. If $k \geq 5$, then $X_\gamma(C_{3k - 1}) = 3k-1 < 2^{k-1} = 2^{\gamma(C_{3k - 1})-1}$, and by Proposition~\ref{prop:gmb>gamma}, we conclude that $\gmb(C_{3k - 1}) > k = \gamma(C_{3k - 1})$. Actually, $\gmb(C_{3k-1})$ is much bigger than $\gamma_g(C_{3k-1})$ as we will see in Section~\ref{sec:cycles}.

The converse of Proposition~\ref{prop:gmb>gamma} does not hold as the following example shows. If $k \in \{3, 4\}$, then $X_{\gamma}(C_{3k-1}) = 3k - 1 > 2^{k-1} = 2^{\gamma(C_{3k - 1})-1}$, but as we will see in Section~\ref{sec:cycles}, $\gamma(C_{3k-1}) = k < k+1 = \left \lfloor \frac{3k-1}{2} \right \rfloor = \gmb(C_{3k-1})$.

\section{Residual graphs}
\label{sec:residual}

In this section we study the Maker-Breaker domination number on a construction that might be of independent interest and that will be later used to determine the invariant for trees. 

If $G$ is a graph, then we say that the {\em residual graph} $R(G)$ of $G$ is the graph obtained from $G$ by iteratively removing pendant paths $P_2$ until no such path is present. By a {\em pendant $P_2$} we mean $P_2$ attached to $G$ with an edge. Hence, when such a pendant $P_2$ is removed, exactly two vertices and two edges are removed. When $G = P_2$, we can also remove it and obtain the empty graph.

Note that $H = R(G)$ for some graph $G$ if and only if $H$ is the empty graph, $H=K_1$, or each support vertex of $H$ has degree at least $3$. This is in particular true if $H$ has no support vertices. We further observe: 

\begin{lemma}
\label{lem:R(T)-is-unique}
If $G$ is a graph, then $R(G)$ is unique (up to isomorphism). 
\end{lemma}

\proof
Let $P'=xy$ be an arbitrary pendant $P_2$ of $G$, where $x$ is a leaf of $G$ and let $z$ be the other neighbor of $y$. Then either $P'$ is removed at some point when pendant $P_2$s are removed from $G$, or $P'$ is not removed at all. The latter possibility can only happen if after this removal process only the path induced with $x,y,z$ remains, and then the pendant path $P''=yz$ is removed. In this case, however, we have $R(G) = K_1$. In addition, by induction every pendant $P_2$ that appears during the removal process will either be eventually removed or will lead to the residual tree $K_1$. 
\qed

Note that the proof of Lemma~\ref{lem:R(T)-is-unique} also reveals that if $R(G)\ne K_1$, then $G\setminus V(R(G))$ is unique. To see that it is not unique in general, consider a path $P_{2k+1}$, $k\ge 2$, and different sequences of removing pendant $P_2$s. 

\begin{lemma}
	\label{lem:R(G)-perfect-matching}
Let $G$ be a graph and $R(G)$ a residual graph of $G$. Then 
\begin{enumerate}
\item[(i)] $G\setminus V(R(G))$ is a forest that has a unique perfect matching, and  
\item[(ii)] $G$ has a perfect matching if and only if $R(G)$ has a perfect matching. 
\end{enumerate}
\end{lemma}

\proof
(i) $G\setminus V(R(G))$ is a forest since it is built from the empty graph by successively attaching to it pendant $P_2$s. If $x_iy_i$, ${i \in I}$, are the pendant $P_2$s that were removed when $R(G)$ was obtained from $G$, then $\{ x_i y_i \}_{i \in I}$ is a unique perfect matching of $G - V(R(G))$. 

(ii) If $G$ has a perfect matching, then its restriction to $G\setminus V(R(G))$ must be the unique perfect matching of $G\setminus V(R(G))$, hence $R(G)$ has a perfect matching. Conversely, if $R(G)$ has a perfect matching, then it can be extended to a perfect matching of $G$ by means of (i). 
\qed

For the proof of the main result of this section, we also need the following. 

\begin{lemma}
	\label{lem:staller}
	If $T$ is a tree that admits a perfect matching and $v \in V(T)$, then Staller has a strategy for the S-game such that Dominator has to select at least $\frac{n(T)}{2}$ vertices to dominate $T$ and $v$ is played by Staller in her last move.
\end{lemma}

\proof
We prove the claim by induction on $n(T)$. If $T=P_2$ and $v \in V(P_2)$, then Staller can play on $v$ and Dominator has to reply on the other vertex.

Let now $n(T)\ge 4$ and consider $T$ as a BFS-tree rooted at an arbitrary  vertex $r$. Let $x$ be a leaf of this BFS-tree at the largest distance from  $r$ and let $y$ be the neighbor of $x$. Then ${\rm deg}(y) = 2$ because $T$ has a perfect matching. Let $z$ be the other neighbor of $y$. Set $T' = T\setminus \{x,y\}$. As $T$ has a perfect matching, $xy$ belongs to it, hence $T'$ also has a perfect matching. If $v \in V(T')$, then Staller starts on $y$, Dominator has to reply on $x$ (otherwise Staller would win) and then Staller applies her strategy on $T'$ (by the induction hypothesis). If $v \in \{x, y\}$, then she applies her strategy on $T'$ with her last move on $z$, and then plays $v$ in her last move. Note that if Dominator plays on $v$ while Staller is playing on $T'$, then Staller wins the game as she can prevent Dominator from playing on one pair of vertices from the matching in $T'$.

From the above strategy of Staller we conclude that the total number of Dominator's moves was $\frac{n(T')}{2} + 1 = \frac{n(T)}{2}$. 
\qed

Note that by the strategy from the proof of lemma~\ref{lem:staller}, unless Staller wants to play on a leaf, she plays on the support vertex, forcing Dominator to reply on its neighboring leaf and separating this $P_2$ from the rest of the graph.

\begin{theorem}
\label{thm:R(G)}
Let $R(G)$ be a residual graph of $G$ and let $H=G\setminus V(R(G))$. Then
\begin{enumerate}
\item[(i)] $\gmb'(G) = \frac{n(H)}{2} + \gmb'(R(G))$,
\item[(ii)]  $\frac{n(H)}{2} + \gmb(R(G)) - 1 \leq \gmb(G) \leq  \frac{n(H)}{2} + \gmb(R(G))$.
\end{enumerate}
\end{theorem}

\proof
(i) $H$ has a perfect matching and is a forest by Lemma~\ref{lem:R(G)-perfect-matching}(i). Let the S-game be played on $G$ and consider the following strategy of Staller. By Lemma~\ref{lem:staller} she can play on each tree of $H$ and play last on the vertex of this tree adjacent to $R(G)$. Dominator has to reply on the matching (otherwise Staller wins the game). Thus, Dominator makes (at least) $\frac{n(H)}{2}$ moves on $H$. Moreover, Staller plays on vertices adjacent to $R(G)$, hence no vertex in $R(G)$ will be dominated by the time Staller makes her first move in $R(G)$. Next, Staller is the player to make the first move on $R(G)$ and she follows her optimal strategy there to  ensure at least $\gmb'(R(G))$ moves of Dominator. 

On the other hand, Dominator's strategy is to then reply wherever Staller plays, $H$ or $R(G)$, with its strategy on this graph. As $H$ has a perfect matching, Dominator makes no more than $\frac{n(H)}{2}$ moves on $H$. Moreover, he makes at most $\gmb'(R(G))$ moves on $R(G)$. Hence, we have $\gmb'(G) = \frac{n(H)}{2} + \gmb'(R(G))$.

(ii) Suppose now that the D-game is played on $G$. To prove the upper bound, Dominator's strategy is to start on $R(G)$ and then reply on $R(G)$ or $H$ if Staller plays there. As $H$ has a perfect matching, Dominator makes no more than $\frac{n(H)}{2}$ moves on $H$. Moreover, he makes at most $\gmb(R(G))$ moves on $R(G)$. Hence we get the upper bound $\gmb(G) \leq \frac{n(H)}{2} + \gmb(R(G))$.

\begin{figure}[ht!]
\centering

\scalebox{0.8}{
\begin{tikzpicture}
\node at (0,0){
\begin{tikzpicture}
\draw (0,0) circle (1.3);
\node at (0,0) {$R(G)$};
\node[noeud, label=left:$d_1$] at (0.7,0.8){};

\draw (0,0) circle (2.5);
\node at (60:2){$H$};

\node[noeud] (a) at (180:1.3){};
\draw (a) -- (165:2.4) -- (195:2.4) -- (a);
\node at (180:2){$T_1$};

\node[noeud] (b) at (230:1.3){};
\draw (b) -- (215:2.4) -- (245:2.4) -- (b);
\node at (230:2){$T_2$};

\node at (280:1.9){\scriptsize $\bullet$};
\node at (290:1.9){\scriptsize $\bullet$};
\node at (300:1.9){\scriptsize $\bullet$};

\node[noeud] (c) at (-10:1.3){};
\draw (c) -- (5:2.4) -- (-25:2.4) -- (c);
\node at (-10:2){$T_k$};

\node at (30:2.8){$G$};

\node at (0,-3){Case 1};

\end{tikzpicture}
};

\node at (7.5,0){
\begin{tikzpicture}
\draw (0,0) circle (1.3);
\node at (0,0) {$R(G)$};

\node[noeud] (a) at (1.3,0){};
\node[noeud] (b) at (2,0){};
\node[noeud] (c) at (3.5,0){};
\node[noeud] (d) at (4.2,0){};
\node[noeud] (e) at (4.9,0){};
\node[noeud] (f) at (5.6,0){};
\node[noeud,label=above:$d_1$] (g) at (6.3,0){};

\node[noeud] (b') at (2,0.7){};
\node[noeud] (d') at (4.2,0.7){};
\node[noeud] (e') at (4.9,0.7){};

\draw (a) -- (b) -- (2.5,0);
\draw[dotted] (2.5,0) -- (3,0);
\draw (3,0) -- (c) -- (d) -- (e) -- (f) -- (g);
\draw (b) -- (b');
\draw (d) -- (d');
\draw (e) -- (e');

\node at (2.65,-2){Case 2.1};

\end{tikzpicture}
};

\node at (3.75,-5){
\begin{tikzpicture}
\draw (0,0) circle (1.3);
\node at (0,0) {$R(G)$};

\node[noeud] (a) at (1.3,0){};
\node[noeud] (b) at (2,0){};
\node[noeud,label=above:$y_2$] (c) at (3.5,0){};
\node[noeud,label=above:$y_1$] (d) at (4.2,0){};
\node[noeud,label=above:$v_{}$] (e) at (4.9,0){};
\node[noeud,label=above:$x_1$] (f) at (5.6,0){};
\node[noeud,label=above:$x_2$] (g) at (6.3,0){};
\node[noeud,label=above:$d_1$] (h) at (7,0){};

\draw (a) -- (b) -- (2.5,0);
\draw[dotted] (2.5,0) -- (3,0);
\draw (3,0) -- (c) -- (d) -- (e) -- (f) -- (g) -- (h);

\node at (2.85,-2){Case 2.2};

\end{tikzpicture}
};

\end{tikzpicture}
}
\caption{Representations of the cases from the proof of Theorem~\ref{thm:R(G)}}
\label{fig:R(G)}
\end{figure}
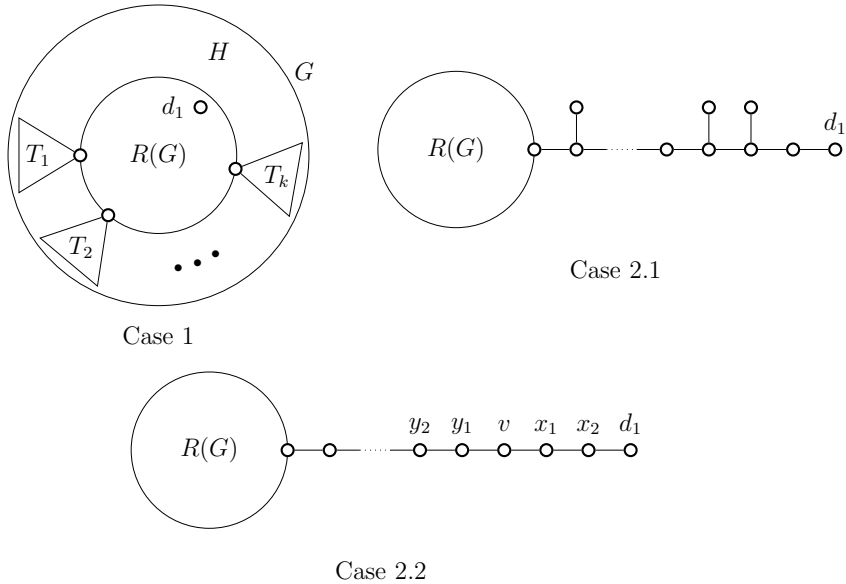

To prove the lower bound, consider the following strategy of Staller depending on the first move of Dominator. We will distinguish two cases, the second with two subcases, which are schematically depicted in Fig.~\ref{fig:R(G)}.  

\begin{description}
	\litem{1} The first move of Dominator is on $R(G)$.
		
	Staller first applies her strategy from Lemma~\ref{lem:staller} on each tree of $H$, playing the vertex adjacent to $R(G)$ as her last move on each of the trees. With this, she forces Dominator to play (at least) $\frac{n(H)}{2}$ moves on $H$. After that we have an ordinary D-game played on $R(G)$, so at least $\gmb(R(G))$ moves are made on it by Dominator if Staller follows her strategy there.
	
	\litem{2} The first move of Dominator is on $H$.
	
	Let $d_1$ be the vertex Dominator plays in his first move, let $T$ be the connected component of $H$ containing $d_1$ (recall that $T$ is a tree), let $P$ be the shortest path between $d_1$ and $R(G)$ in $T$, and let $M$ be the unique perfect matching of $T$ (cf.~Lemma~\ref{lem:R(G)-perfect-matching}(i)).
	
	In this case, Staller first applies her strategy from Lemma~\ref{lem:staller} on all the other trees of $H$, playing the vertex adjacent to $R(G)$ as her last move on each tree. Next, Staller applies her strategy from Lemma~\ref{lem:staller} on the edges from $M$, which are not incident with $P$. Additionally, she plays last on the vertices closest to $P$. After that, only $R(G)$, $P$, and maybe some vertices adjacent to $P$, remain undominated. 
	
	\begin{description}
		\litem{2.1} At least one vertex adjacent to $P$ is still undominated (see Fig.~\ref{fig:R(G)}). 
		
		Let $u$ be an undominated vertex adjacent to $P$. Staller plays on its neighbor on $P$, forcing Dominator to reply on $u$. Staller does so on each such vertex. After that, the only undominated vertices lie on $P$, moreover, up to now, at least one move of Dominator was played on each already completely dominated edge from $M$. 
		
		As long as there are some more undominated edges from $M$ on $P$, at least one of them, say $e \in M$,  is adjacent to a vertex $s$ of $P$ already played by Staller. Her strategy is to play on the vertex of $e$ which is at distance $2$ from $s$. Then Dominator has to reply on the other vertex of $e$, otherwise Staller wins by playing it. Hence, Staller can force Dominator to reply on all remaining edges.
		
		\litem{2.2} The only undominated vertices in $H$ lie on $P$. 
		
		Staller's strategy is to play on the vertex $v$ of $P$ at distance $3$ from $d_1$. Dominator has to reply on a neighbor of $v$, otherwise one of the neighbors of $v$ is not dominated, say $u$, and Staller can win by playing $u$ and later playing the unplayed vertex from $N[v]$ or $N[u]$. Indeed, in this case, two undominated adjacent vertices are played by Staller, and no matter where Dominator answers, she can play another consecutive vertex and win the game. 
		
		Let $x_i$ be a vertex at distance $i$ from $v$ on $P$ in the direction of $d_1$, and $y_i$ be a vertex at distance $i$ from $v$ on $P$ in the direction of $R(G)$ for all possible $i \geq 1$, see Fig.~\ref{fig:R(G)} again. 
		
		If Dominator replies on $x_1$, then Staller's next move is $y_2$. Now, Dominator has to reply on $y_1$, otherwise Staller wins. Then Staller repeats this strategy until $P$ is dominated, i.e., she plays on the vertices $y_{2k}$ in the increasing order, and Dominator is forced to reply on $y_{2k-1}$.
		
		If Dominator replies on $y_1$, then Staller replies on $x_2$. After that, Dominator has to play $x_1$. Next, Staller applies the same strategy as before, taking $y_1$ as the new $d_1$. 		
	\end{description}
	
	In both cases, Dominator is forced to play at least one move on each edge of the matching $M$, hence at least $\frac{n(T)}{2}$ moves are made on $T$. On $H-T$, at least $\frac{n(H - T)}{2}$ moves are made by Lemma~\ref{lem:staller}.

	After $T$ is completely dominated, Staller follows her optimal strategy on $R(G)$, but it might happen that one vertex $u$ in $R(G)$ is already dominated (by a move of Dominator in $H$ close to $R(G)$). As Staller's strategy on $H$ forces Dominator to answer on $H$, Staller will be the first player to play on $R(G)$. But as she can imagine that Dominator's move was $u$, we have 
	$$\gmb'(R(G)|u) \geq \gmb'(R(G)|N[u]) \geq \gmb(R(G))-1\, ,$$
	hence the total number of moves on $R(G)$ is at least $\gmb(R(G)) - 1$. 	
\end{description}

In either case, Dominator played at least $\frac{n(H)}{2} + \gmb(R(G)) - 1$ moves, which proves the lower bound.
\qed

Note that in the inequality $\gmb'(G|u) \geq \gmb(G)-1$ from the above proof, equality can be attained. For example, consider the graph $G$ on Fig.~\ref{fig:lowerR(G)}. Clearly, $\gmb(G) = 2$ and $\gmb'(G|u) = 1$.

\begin{figure}[ht]
	\begin{center}
		\begin{tikzpicture}
		
		\node[noeud] (0) at (0,0){};
		\node[noeud] (1) at (0.866, 0.5){};
		\node[noeud] (2) at (0.866, -0.5){};
		\node[noeud] (3) at (1.732, 0){};
		\node[noeud, label=right:$u$] (4) at (2.732, 0){};
				
		\draw (0) -- (1) -- (3) -- (4);
		\draw (0) -- (2) -- (3);
		\draw (1) -- (2);

		\end{tikzpicture}
	\end{center}
	\caption{A graph $G$ with the property $\gmb'(G|u) = \gmb(G)-1$.}
	\label{fig:lowerR(G)}
\end{figure}

Appending to $G$ some trees with perfect matchings, where at least one of them is attached to $u$, we get graphs that attain the lower bound from Theorem~\ref{thm:R(G)}(ii).

To conclude the section we apply the residual construction to determine the Maker-Breaker domination number of trees. This contrasts the domination game where no such result is known, cf.~\cite{bresar-2013, henning-2017, henning-2017b}.

\begin{theorem}
	\label{thm:trees}
	If $T$ a tree, then 
	$$\gmb(T) = \begin{cases}
	\frac{n(T)}{2}; & T \text{ has a perfect matching},\\
	\frac{n(T)-1}{2}; & R(T) \cong K_1,\\
	\frac{n(T)-k+1}{2}; & R(T) \cong K_{1,k} \text{ for } k \geq 3,\\
	\infty; & \text{otherwise},
	\end{cases}$$
	and 
	$$\gmb'(T) = \begin{cases}
	\frac{n(T)}{2}; & T \text{ has a perfect matching},\\
	\infty; & \text{otherwise}.
	\end{cases}$$
\end{theorem}

\proof
Let $T$ be a tree. Note first that by Lemma~\ref{lem:R(G)-perfect-matching} the cases in Theorem~\ref{thm:trees} are disjoint.

If $T$ admits a perfect matching, then $\gmb(T) \leq \frac{n(T)}{2}$ since by playing once on every edge of the matching, Dominator dominates $T$. By Lemma~\ref{lem:staller}, $\gmb(T) \geq \frac{n(T)}{2}$. Hence, $\gmb(T) = \frac{n(T)}{2}$ and by the same reasoning, $\gmb'(T) = \frac{n(T)}{2}$.

If $R(T) \cong K_1$ and $v$ is the remaining vertex, then $T - v$ admits a perfect matching by Lemma~\ref{lem:staller}. In this case, Dominator first plays on a neighbor of $v$ (which also belongs to one edge of the matching in $T-v$) and then follows the perfect matching in $T-v$. Thus he makes at most $\frac{n(T)-1}{2}$ moves. By Lemma~\ref{lem:staller}, Staller can force him to play on the perfect matching except for his first move, hence $\gmb(T) = \frac{n(T) - 1}{2}$. 

If $R(T) \cong K_{1,k}$ for some $k \geq 3$, then $T - K_{1,k}$ admits a perfect matching. Dominator can start the game on the center vertex of the star $K_{1,k}$ and then follow the matching in the remaining part. Thus he makes at most $1 + \frac{n(T)-k-1}{2} = \frac{n(T)-k+1}{2}$ moves. On the other hand, Staller can ensure that he cannot win in less moves. If his first move is indeed in the center of a star, then Lemma~\ref{lem:staller} assures that at least $1 + \frac{n(T)-k-1}{2}$ Dominator's moves are needed. If his first move is elsewhere, then Staller can ensure that Dominator follows her moves in the perfect matchings in subtrees of $T - K_{1,k}$ while additionally she plays the last move on a vertex closest to $K_{1,k}$ in all subtrees of at least two descendants of the star (again by Lemma~\ref{lem:staller}). With this she forces Dominator to play on the center of the star and Dominator makes at least $1 + \frac{n(T)-k-1}{2}$ moves. If Dominator played elsewhere, then Staller would win by playing all vertices in a closed neighborhood of one of those two descendants of the star. This proves that $\gmb(T) = \frac{n(T)-k+1}{2}$.

Otherwise, we know by~\cite{gledel-2018+} that Staller wins on $T$ (no matter which player starts the game), hence $\gmb(T) = \gmb'(T) = \infty$.
\qed

Note that by Theorem~\ref{thm:trees}, $\gmb$ and $\gmb'$ of trees can be computed in polynomial time. 

\section{Cycles}
\label{sec:cycles}


The D-game domination number and the S-game domination number of cycles are given with the following formulas:
$$\gamma_g(C_n) = \begin{cases}
\left \lceil \frac{n}{2} \right \rceil - 1; & n \equiv 3 \mod 4,\\
\left \lceil \frac{n}{2} \right \rceil; & \text{otherwise},
\end{cases}
\qquad
\gamma_g'(C_n) = \begin{cases}
\left \lceil \frac{n-1}{2} \right \rceil - 1; & n \equiv 2 \mod 4,\\
\left \lceil \frac{n-1}{2} \right \rceil; & \text{otherwise}.
\end{cases}
$$
This fundamental result was first obtained in an unpublished manuscript~\cite{kinnersley-2012}. The result appeared for the first time in press in the paper~\cite{kosmrlj-2017}, where an alternative proof is given. For the total domination game, parallel results were obtain in~\cite{dorbec-2016}. The latter paper investigates the total domination game on paths and cycles only. So the (total) game domination number of cycles is far from being straightforward. Here we determine the Maker-Breaker domination number of cycles, a task that turned out to be less involved. 

\begin{theorem}
\label{thm:cycles}
If $n \geq 3$, then $$\gmb(C_n) = \gmb'(C_n) = \left \lfloor \frac{n}{2} \right \rfloor\,.$$
\end{theorem}

\proof
We begin by showing that $\gmb(C_n) \leq \left \lfloor \frac{n}{2} \right \rfloor$ and $\gmb'(C_n) \leq \left \lfloor \frac{n}{2} \right \rfloor$. If $n$ is even, $C_n$ has a perfect matching (which is also a dominating matching), thus by Fact~\ref{fact:dominating-pair}, it holds that $\gmb(C_n) \leq \frac{n}{2} = \left \lfloor \frac{n}{2} \right \rfloor$ and $\gmb'(C_n) \leq \left \lfloor \frac{n}{2} \right \rfloor$. Now consider the case when $n$ is odd. In a D-game, let $v$ be the first vertex played by \Dom. Clearly, among undominated vertices, $V(C_n) - N[v]$, there is a perfect matching. Thus $\gmb(C_n) \leq 1 + \frac{n-3}{2} = \left \lfloor \frac{n}{2} \right \rfloor$. In an S-game, suppose $s_1' = u$. Then \Dom\ should reply on a neighbor $v$ of the vertex $u$. Now there is a perfect matching among $V(C_n) - N[v]$, so $\gmb'(C_n) \leq 1 + \frac{n-3}{2} = \left \lfloor \frac{n}{2} \right \rfloor$. This proves the upper bounds. 

To find the lower bounds we need to find an appropriate strategy for \St. Set for the rest of the proof that $V(C_n) = \{x_1, \ldots, x_n\}$, where the adjacencies are natural. 

We first show the lower bound for the S-game: $\gmb'(C_n) \geq \left \lfloor \frac{n}{2} \right \rfloor$. Suppose, without loss of generality, that $s_1' = x_1$. Notice that \Dom\ has to reply on a neighbor of $x_1$, for otherwise \St\ plays as $s_2'$ the not yet dominated neighbor of $x_1$. Then \Dom\ cannot in one move dominate $s_1'$ and $s_2'$. Say he leaves $s_2'$ undominated. Then \St\ can play the other neighbor of $s_2'$ and win the game as \Dom\ cannot dominate $s_2'$. So without loss of generality, \Dom\ replies with $d_1' = x_n$. \St's next move is $s_2' = x_3$. In order to prevent \St\ from winning, \Dom\ has to play $d_2' = x_2$. Then \St\ continues with the same strategy, forcing \Dom\ to play on (almost) every second move. 

If $n$ is even, the game ends after \St\ plays $x_{n-1}$ and \Dom\ replies on $x_{n-2}$. So in this case, \Dom\ plays all even labeled vertices, hence $\gmb'(C_n) \geq \frac{n}{2} = \left \lfloor \frac{n}{2} \right \rfloor$. If $n$ is odd, the game ends after \St\ plays $x_{n-2}$ and \Dom\ replies on $x_{n-3}$ (as $x_{n-1}$ is already dominated by $d_1'$). So \Dom\ again plays all even labeled vertices, thus $\gmb'(C_n) \geq \frac{n-1}{2} =  \left \lfloor \frac{n}{2} \right \rfloor$. 

It remains to prove that $\gmb(C_n) \geq \left \lfloor \frac{n}{2} \right \rfloor$. Assume without loss of generality that $d_1 = x_{1}$. \St\ replies on $s_1 = x_4$ (so at distance $3$ from $d_1$). Again, \Dom\ has to reply on a neighbor of $s_1$. If he replies on $x_3$, then \St\ can apply the above strategy by playing $x_6$ next, and then every second vertex along the cycle to ensure at least $\left \lfloor \frac{n}{2} \right \rfloor$ moves. 

If \Dom\ replies on $d_2 = x_5$, then \St\ plays $s_2 = x_8$ (again at distance $3$ from $d_2$ in the same direction as before). But observe that at some point of the game \Dom\ will have to play on $\{x_2, x_3\}$ to dominate the whole graph. 

By repeating this strategy, \St\ ensures that among every four consecutive vertices of the cycle, \Dom\ makes at least two moves (except maybe in the last one, two or three remaining vertices). We now distinguish four different cases. 
\begin{itemize}
	\item If $n \equiv 0 \bmod 4$, then no vertex remains and $\gmb(C_n) \geq \frac{2 n}{4} = \left \lfloor \frac{n}{2} \right \rfloor$.
	\item If $n \equiv 1 \bmod 4$, then only one vertex remains, which is already dominated by $d_1$, so $\gmb(C_n) \geq \frac{2 (n-1)}{4} = \left \lfloor \frac{n}{2} \right \rfloor$.
	\item If $n \equiv 2 \bmod 4$, then among the remaining two vertices, one is dominated by $d_1$ but the other is not. So \Dom\ has to make another move, thus $\gmb(C_n) \geq \frac{2(n-2)}{4} + 1 = \left \lfloor \frac{n}{2} \right \rfloor$.
	\item If $n \equiv 3 \bmod 4$, then two of the remaining three vertices are not yet dominated, so \Dom\ still has to make just one more move. So $\gmb(C_n) \geq \frac{2(n-3)}{4} = \left \lfloor \frac{n}{2} \right \rfloor$. 
\qed
\end{itemize}

\section{Union of graphs}
\label{sec:union}

The Maker-Breaker domination game was in~\cite{gledel-2018+} studied on disjoint unions on graphs, the obtained results were in particular applied to cographs. In this section we complement their investigation with the following result. 

\begin{theorem}
\label{thm:union}
If $G$ and $H$ are graphs, then
$$(i)\ \gmb(G) + \gmb(H) \leq \gmb(G \cup H) \leq \min\{\gmb'(G) + \gmb(H),\gmb(G) + \gmb'(H)\}\,,$$
$$(ii)\ \max\{\gmb'(G) + \gmb(H), \gmb(G) + \gmb'(H)\} \leq \gmb'(G \cup H) \leq \gmb'(G) + \gmb'(H)\,.$$
Moreover, all the bounds are sharp.
\end{theorem}

\proof
(i) It follows directly from the results from~\cite{gledel-2018+} that if $\gmb(G) = \infty$ or $\gmb(H) = \infty$, then $\gmb(G\cup H) = \infty$. Suppose then that both $\gmb(G)$ and $\gmb(H)$ are finite. We give a strategy for Staller such that when a D-game is played on $G \cup H$, she can ensure that Dominator will select at least $\gmb(G) + \gmb(H)$ vertices. The strategy of Staller is the following: each time Dominator selects a vertex from $G$ or from $H$, she answers optimally in $G$ or $H$ (with an optimal strategy restricted to $G$ or $H$), respectively, as long as this is possible. Suppose without loss of generality that Dominator has first dominated $G$. Then Staller can reply either in $G$, provided she has a legal move in $G$ available, or in $H$. In the first case a usual D-game will be played on $H$. In the second case we have a game on $H$ in which Dominator has passed one move. By the No-Skip Lemma, Dominator will need to select at least $\gmb(H)$ vertices from $H$. In any case, $\gmb(G \cup H) \geq \gmb(G) + \gmb(H)$.

To prove the upper bound, suppose first that Dominator starts by playing his optimal move on $G$ and then follows Staller in $G$ or in $H$ whenever she plays in $G$ or in $H$. In this way (having in mind the No-Skip Lemma) Dominator achieves at most $\gmb(G)$ moves in $G$ and at most $\gmb'(H)$ moves in $H$, hence at most $\gmb(G)+\gmb'(H)$ vertices for $G \cup H$. If instead he starts by playing his optimal move on $H$, then he can guarantee to play at most $\gmb'(G)+\gmb(H)$ vertices. By choosing the smaller of the two values, Dominator has a strategy such that he dominates $G \cup H$ with no more than $\min\{\gmb'(G) + \gmb(H),\gmb(G) + \gmb'(H)\}$ moves and so $\gmb(G \cup H) \leq \min\{\gmb'(G) + \gmb(H),\gmb(G) + \gmb'(H)\}$.

(ii) This is done using similar arguments as in (i).

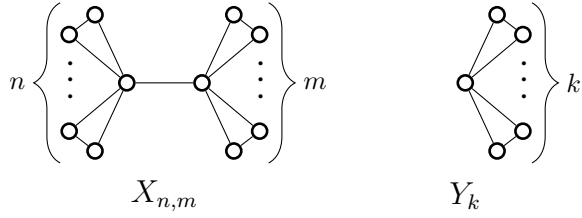
\begin{figure}[ht]
\begin{center}
\begin{tikzpicture}

\node at (0.5,-1.5){$X_{n,m}$};

\node[noeud] (x1) at (0,0){};

\node[noeud] (a1) at (115:1){};
\node[noeud] (b1) at (140:1){};
\node[noeud] (an) at (220:1){};
\node[noeud] (bn) at (245:1){};

\node[noeud] (x2) at (1,0){};
\node[noeud] (u1) at ($ (1,0)+(65:1)$){};
\node[noeud] (v1) at ($ (40:1)+(1,0) $){};
\node[noeud] (um) at ($ (-40:1)+(1,0) $){};
\node[noeud] (vm) at ($ (-65:1)+(1,0) $){};

\draw (x1) -- (x2);
\draw (x1) -- (a1);
\draw (x1) -- (an);
\draw (x1) -- (b1);
\draw (x1) -- (bn);
\draw (a1) -- (b1);
\draw (an) -- (bn);
\draw (x2) -- (u1);
\draw (x2) -- (um);
\draw (x2) -- (v1);
\draw (x2) -- (vm);
\draw (u1) -- (v1);
\draw (um) -- (vm);

\path (b1) -- (an) node[midway, sloped] {\large $\cdots$};
\path (v1) -- (um) node[midway, sloped] {\large $\cdots$};

\draw [decorate,decoration={brace,amplitude=10pt,raise=4pt}] (235:1.3) -- (125:1.3) node[left= 13pt,midway]{\footnotesize $n$};
\draw [decorate,decoration={brace,amplitude=10pt,raise=4pt}] ($ (1,0)+(55:1.3)$) -- ($ (-55:1.3)+(1,0) $) node[right= 13pt,midway]{\footnotesize $m$};

\begin{scope}[xshift = 4.5cm]

\node at (0,-1.5) {$Y_k$};

\node[noeud] (y) at (0,0){};
\node[noeud] (c1) at (65:1){};
\node[noeud] (d1) at (40:1){};
\node[noeud] (ck) at (-40:1){};
\node[noeud] (dk) at (-65:1){};

\draw (y) -- (c1){};
\draw (y) -- (ck){};
\draw (y) -- (d1){};
\draw (y) -- (dk){};

\draw (c1) -- (d1);
\draw (ck) -- (dk);
\path (d1) -- (ck) node[midway, sloped] {\large $\cdots$};

\draw [decorate,decoration={brace,amplitude=10pt,raise=4pt}] (55:1.3) -- (-55:1.3) node[right= 13pt,midway]{\footnotesize $k$};

\end{scope}

\end{tikzpicture}
\end{center}
\caption{Representation of graphs $X_{n,m}$ and $Y_k$}
\label{fig:XnmYk}
\end{figure}

To demonstrate the sharpness of the bounds, consider the graphs $X_{n,m}$, $1 \leq m \leq n$, and $Y_k$, $k\ge 1$, as depicted in Fig.~\ref{fig:XnmYk}.  

Observe first that $\gmb(X_{n,m}) = m+1$, $\gmb'(X_{n,m}) = n+1$, $\gmb(Y_{k}) = 1$, and $\gmb'(Y_{k}) = k$. Consider next the union $X_{n,m} \cup Y_{k}$, where $k\le m$, and the union $Y_{k_1} \cup Y_{k_2}$, where $k_1 \leq k_2$. Then in the D-game as well in the S-game played in either of the unions, an optimal first move is to play a vertex of highest degree. Moreover, an optimal reply to this move is to play a vertex of second highest degree. Therefore, 
\begin{itemize}
\item $\gmb(X_{n,m} \cup Y_{k}) = 1 + m +1$, reaching the lower bound of (i), 
\item $\gmb'(X_{n,m} \cup Y_{k}) = n + 1 +k $, reaching the upper bound of (ii),
\item $\gmb(Y_{k_1} \cup Y_{k_2}) = k_1 +1$, reaching the upper bound of (i), and 
\item $\gmb'(Y_{k_1} \cup Y_{k_2}) = 1 + k_2$, reaching the lower bound of (ii).
\qed
\end{itemize}

\section{Concluding remarks}
\label{sec:concluding}

To conclude the paper we list several problems and directions for further investigation of the Maker-Breaker domination number. 

\begin{enumerate}
\item For the upper bound in~\eqref{eq:n/2-Dominator} we have provided examples of graphs that attain equality. These examples are not connected and it is  not difficult to achieve equality with connected graphs of even order. However, we do not know of any connected graph of odd order (different from $K_1$) for which equality in~\eqref{eq:n/2-Dominator} is achieved. More generally, we ask for a characterization of the extremal graphs with respect to~\eqref{eq:n/2-Dominator} and~\eqref{eq:n/2-Staller}. 

\item As we already mentioned, it would be interesting to find a structural characterization of the graphs $G$ for which $\gmb(G) = \gamma(G) = k$ holds, where $k\ge 2$ is a fixed integer. 

\item It would also be interesting to investigate $\gmb(G\cp H)$ and $\gmb'(G\cp H)$, where $G$ and $H$ are arbitrary graphs and $G\cp H$ is the Cartesian product of $G$ and $H$. In particular, it would be interesting to determine $\gmb(P_n\cp P_m)$ (and $\gmb'(P_n\cp P_m)$), as well as $\gmb(G\cp K_2)$ (and $\gmb'(G\cp K_2)$) for an arbitrary graph $G$. 

\item If $G$ is a cograph, then it is not difficult to determine whether Dominator or Staller wins the Maker-Breaker domination game~\cite{gledel-2018+}. On the other hand, it does not seem straightforward to determine the Maker-Breaker domination numbers of co-graphs. 

\item In this paper we have considered the Maker-Breaker domination number which is an optimization problem from Dominator's point of view. It would likewise be of interest to consider the Staller's point of view, that is, assuming that Staller wins on a graph $G$, what is the minimum number of moves with which she can achieve the goal? 

\end{enumerate}

\section*{Acknowledgements}

We acknowledge the financial support from the Slovenian Research Agency (research core funding No.\ P1-0297 and projects J1-9109, N1-0095).	


\end{document}